\documentclass{gtart}

\def\ifplaintex{\expandafter\ifx\csname documentclass\endcsname\relax}

\def\gtp{{\mathsurround=0pt\it $\cal G\mskip-2mu$eometry \&\ 
$\cal T\!\!$opology $\cal P\!$ublications}}  

\def\recd{{\small Received:\qua\receiveddate\ifx\reviseddate\relax
\else\qquad Revised:\qua\reviseddate\fi\par}} 

\def\lognumber#1{\def\thelognumber{#1}}
\def\volumenumber#1{\def\thevolumenumber{#1}}
\def\volumeyear#1{\def\thevolumeyear{#1}}
\def\papernumber#1{\def\thepapernumber{#1}}
\def\pagenumbers#1#2{\def\startpage{#1}\def\finishpage{#2}}
\def\published#1{\def\publishdate{#1}}

\def\received#1{\def\receiveddate{#1}}
\def\revised#1{\def\reviseddate{#1}}
\def\accepted#1{\def\accepteddate{#1}}

\def\asciiauthors#1{\def\theasciiauthors{#1}}
\def\asciiaddress#1{\def\theasciiaddress{#1}}

\def\coverauthors#1{\def\thecoverauthors{#1}}
\long\def\asciiabstract#1{\long\def\theasciiabstract{#1}}

\let\\\par\let\thelognumber\relax\let\thevolumenumber\relax
\let\thepapernumber\relax\let\thevolumeyear\relax\let\startpage\relax
\let\finishpage\relax\let\publishdate\relax\let\receiveddate\relax
\let\reviseddate\relax\let\accepteddate\relax\let\theasciititle\relax
\let\theasciiauthors\relax\let\theasciiaddress\relax
\let\theasciiabstract\relax

\let\thecoverauthors\relax\let\theasciiemail\relax

\ifplaintex
\font\logobig=cmssbx10 scaled 3836
\font\logomed=cmssbx10 scaled 2557
\else
\font\logobig=cmssbx10 scaled 4200
\font\logomed=cmssbx10 scaled 2800
\fi

\long\def\makeagttitle{   
\count0=\startpage
\agt\hfill      
\hbox to 45truept{\vbox to 0pt{\vglue -13truept{\logomed A\kern -.37em{\logobig 
T}\kern -.38em G}\vss}\hss}
\break
{\small Volume \thevolumenumber\ (\thevolumeyear)
\startpage--\finishpage\nl
Published: \publishdate}

\vglue .25truein

{\parskip=0pt\leftskip 0pt plus
1fil\def\\{\par\smallskip}{\Large\bf\thetitle}\par\medskip} \vglue
0.05truein

{\parskip=0pt\leftskip 0pt plus 1fil\def\\{\par}{\sc\theauthors}
\par\medskip}%
 
\vglue 0.03truein 

{\small\leftskip 25truept\rightskip 25truept{\bf Abstract}\stdspace\theabstract

{\bf AMS Classification}\stdspace\theprimaryclass
\ifx\thesecondaryclass\relax\else; \thesecondaryclass\fi\par
{\bf Keywords}\stdspace \thekeywords\par}\vglue 7truept

}   

\ifplaintex
\font\phead=cmsl9 scaled 950
\font\pnum=cmbx10 scaled 913
\font\pfoot=cmsl9 scaled 950
\headline{\vbox to 0pt{\vskip -4.5mm\line{\small\phead\ifnum
\count0=\startpage ISSN 1472-2739 (on-line) 1472-2747 (printed)
\hfill {\pnum\folio}\else\ifodd\count0\def\\{ }%
\ifx\theshorttitle\relax\thetitle\else\theshorttitle\fi\hfill{\pnum\folio}
\else\def\\{ and }{\pnum\folio}\hfill\ifx\theshortauthors\relax\theauthors
\else\theshortauthors\fi\fi\fi}\vss}}
\footline{\vbox to 0pt{\vglue 0mm\line{\small\pfoot\ifnum\count0=\startpage
\copyright\ \gtp\hfill\else
\agt, Volume \thevolumenumber\ (\thevolumeyear)\hfill\fi}\vss}}
\else
\font=cmsl9 scaled 1050
\font\lnum=cmbx10 
\font=cmsl9 scaled 1050
\makeatletter
\def\@oddhead{{\small\lhead\ifnum\count0=\startpage ISSN 1472-2739 
(on-line) 1472-2747 (printed)\hfill {\lnum\number\count0}\else\ifodd\count0
\def\\{ }\ifx\theshorttitle\relax \thetitle \else\theshorttitle\fi\hfill
{\lnum\number\count0}\else\def\\{ and }{\lnum\number\count0}
\hfill\ifx\theshortauthors\relax 
\theauthors\else\theshortauthors\fi\fi\fi}}\def\@evenhead{\@oddhead}
\def\@oddfoot{\small\lfoot\ifnum\count0=\startpage\copyright\ \gtp\hfill\else
\agt, Volume \thevolumenumber\ (\thevolumeyear)\hfill\fi}
\def\@evenfoot{\@oddfoot}
\makeatother
\fi
\let\maketitlepage\makeagttitle
\let\makeshorttitle\maketitlepage
\let\maketitle\maketitlepage

\newwrite\gtoutfile
\long\gdef\makeheadfile{  
{\def\\{, }\def\s{ }
\immediate\openout\gtoutfile head.xxx
\immediate\write\gtoutfile{To: math@arxiv.org}
\immediate\write\gtoutfile{Subject: put OR rep NNNNN:ppppp}
\immediate\write\gtoutfile{--text follows this line--}
\immediate\write\gtoutfile{Proxy-for: \ifx\theasciiauthors\relax
\theauthors\else\theasciiauthors\fi\s<\ifx\theasciiemail\relax\theemail\else\theasciiemail\fi>}
\immediate\write\gtoutfile{\noexpand\\}
\immediate\write\gtoutfile{Authors: \ifx\theasciiauthors\relax
\theauthors\else\theasciiauthors\fi}
{\def\\{ }\immediate\write\gtoutfile{Title: \ifx\theasciititle\relax
\thetitle\else\theasciititle\fi}}
\immediate\write\gtoutfile{Subj-class: GT or SG, GR etc}
\immediate\write\gtoutfile{MSC-class: \theprimaryclass\ifx\thesecondaryclass\relax\else, \thesecondaryclass\fi}
\immediate\write\gtoutfile{Journal-ref: Algebr. Geom. Topol. \thevolumenumber\s
(\thevolumeyear) \startpage-\finishpage}
\immediate\write\gtoutfile{Comments: Published by Algebraic and
Geometric Topology at}
\immediate\write\gtoutfile{\s\s\s  http://www.maths.warwick.ac.uk/agt/AGTVol\thevolumenumber/agt-\thevolumenumber-\thepapernumber.abs.html}
\immediate\write\gtoutfile{\noexpand\\}
\immediate\write\gtoutfile{}
\ifx\theasciiabstract\relax
\immediate\write\gtoutfile{\theabstract}\else
\immediate\write\gtoutfile{\theasciiabstract}\fi
\immediate\write\gtoutfile{}
\immediate\write\gtoutfile{\noexpand\\}
\immediate\write\gtoutfile{}
\immediate\closeout\gtoutfile}}  

\def\maketitlepage{\makeagttitle\makeheadfile}
\let\makeshorttitle\maketitlepage
\let\maketitle\maketitlepage

\lognumber{39}
\volumenumber{3}
\volumeyear{2003}
\papernumber{39}
\published{1 November 2003}
\pagenumbers{1103}{1118}
\received{29 June 2003}
\revised{16 October 2003}
\accepted{20 October 2003}

\usepackage{amssymb,graphicx}

\theoremstyle{plain}
\newtheorem{theorem}{Theorem}[section]
\newtheorem{lemma}[theorem]{Lemma}
\newtheorem{proposition}[theorem]{Proposition}
\newtheorem{corollary}[theorem]{Corollary}

\theoremstyle{definition}
\newtheorem{definition}[theorem]{Definition}

\theoremstyle{remark}

\def\CRS{\mbox{CRS}}

\begin{document}

\title{The $n$th root of a braid is unique up to conjugacy}
\author{Juan Gonz\'alez-Meneses}
\asciiauthors{Juan Gonzalez-Meneses}
\coverauthors{Juan Gonz\noexpand\'alez-Meneses}
\address{Universidad de Sevilla. Dep. Matem\'atica Aplicada I\\ETS Arquitectura, Av. Reina Mercedes 2, 41012-Sevilla, Spain}
\asciiaddress{Universidad de Sevilla. Dep. Matematica Aplicada I\\ETS Arquitectura, Av. Reina Mercedes 2, 41012-Sevilla, Spain}

\email{meneses@us.es}

\url{www.personal.us.es/meneses}

\begin{abstract}
We prove a conjecture due to Makanin: if $\alpha$ and $\beta$ are elements of the Artin
braid group $B_n$ such that $\alpha^k=\beta^k$ for some nonzero integer $k$, then $\alpha$
and $\beta$ are conjugate. The proof involves the Nielsen-Thurston classification of braids.
\end{abstract}

\asciiabstract{We prove a conjecture due to Makanin: if \alpha and
\beta are elements of the Artin braid group B_n such that
\alpha^k=\beta^k for some nonzero integer k, then \alpha and \beta are
conjugate. The proof involves the Nielsen-Thurston classification of
braids.}



\maketitle

\section{Introduction}

The Artin braid group on $n$ strands, $B_n$, is the group of automorphism of the
$n$-punctured disc that fix the boundary pointwise, up to isotopies relative to the
boundary. One can also consider the elements of $B_n$ ({\em braids}) as isotopy classes of
loops in the space of configurations of $n$ points in a disc $D$. That is, a braid is
represented by the disjoint movements of $n$ points in the disc, starting and ending with
the same configuration, maybe permuting their positions. Braids can also be represented in a
three dimensional picture: if we consider the cylinder $D\times [0,1]$, and fix $n$ base
points $P_1,\ldots,P_n$ in $D$, the movement of the point $P_i$ is represented by a path,
called $i$th {\em strand}, going from $P_i\times \{0\}$ to some $P_j\times \{1\}$. The $n$
strands of a braid are always disjoint, and isotopies correspond to deformations of the
strands, keeping the endpoints fixed.

 The braid group $B_n$ is of interest in several fields of mathematics, with important
applications to low dimensional topology, knot theory, algebraic geometry or
cryptography. Among the basic results concerning braid groups, one can find the
presentation, in terms of generators and relations, given by Artin~\cite{A}
$$
B_{n}=\left\langle \sigma _{1},\sigma _{2},\ldots ,\sigma _{n-1}  \left|
\begin{array}{ll}
\sigma _{i}\sigma _{j}=\sigma _{j}\sigma _{i} & (|i-j|\geq 2) \\
\sigma _{i}\sigma _{i+1}\sigma _{i}=\sigma _{i+1}\sigma _{i}\sigma _{i+1} & (1 \leq i
\leq n-2)
\end{array}
\right. \right\rangle,   \label{presen}
$$
and the solutions to the word problem~\cite{Garside,ThBnaut,BKL} and the conjugacy
problem~\cite{EM,BKL,FrGM_conj}. Closely related to the latter is the problem of determining
the centralizer of a given braid~\cite{Makanin,GMFr,GMW}. It is in the context of these two
problems (conjugacy problem and computation of centralizers) that extraction of roots in
braid groups becomes interesting (see Corollary~\ref{C:main} in this paper, and the last
section of~\cite{GMW}).

\medskip
{\bf Remark}\qua We read the product of two braids $\alpha \beta$ from left to right,
as is usually done in braid theory. That is, if we consider $\alpha$ and $\beta$ as
automorphisms of the disc $D$, then $\alpha \beta (D)=\beta(\alpha(D))$.
\medskip

Given a braid $\eta\in B_n$ and an integer $k$, the problem to determine if there exists
$\alpha\in B_n$ such that $\alpha^k=\eta$ has been solved in~\cite{Sty} (see
also~\cite{Sibert}). But it is known that such an $\alpha$ is not necessarily unique. For
instance $(\sigma_1\sigma_2)^3=(\sigma_2\sigma_1)^3$ but $\sigma_1\sigma_2 \neq
\sigma_2\sigma_1$. G. S. Makanin (see~\cite{NYGTC}, problem B11) conjectured that any two
solutions of the above equation are conjugate, and in this paper we will show that it is
true. In other words, we show the following:

\begin{theorem}\label{T:main}
If $\alpha,\beta\in B_n$ are such that $\alpha^k=\beta^k$ for some $k\neq 0$, then
$\alpha$ and $\beta$ are conjugate.
\end{theorem}

Our proof involves the Nielsen-Thurston classification of braids into periodic, reducible or
pseudo-Anosov. We will see, for instance, that the $k$th root of a pseudo-Anosov braid is
unique, if it exists, while a periodic or reducible braid may have several roots.

One easy consequence of Theorem~\ref{T:main} is the following, which could be useful for
testing conjugacy in braid groups.

\begin{corollary}\label{C:main}
Let $\alpha,\beta\in B_n$ and let $k$ be a nonzero integer. Then $\alpha$ is conjugate to
$\beta$ if and only if $\alpha^k$ is conjugate to $\beta^k$.
\end{corollary}

\begin{proof}
If $\alpha$ is conjugate to $\beta$ then $\eta^{-1}\alpha \eta=\beta$ for some braid $\eta$.
Then $\eta^{-1}\alpha^k \eta = (\eta^{-1}\alpha \eta)^k = \beta^k$, hence $\alpha^k$ and
$\beta^k$ are conjugate.

Conversely, suppose that $\alpha^k$ and $\beta^k$ are conjugate. Then $\eta^{-1}\alpha^k
\eta = \beta^k$ for some $\eta$. This means that $(\eta^{-1}\alpha \eta)^k = \beta^k$, and
by Theorem~\ref{T:main} this implies that $\beta$ is conjugate to $\eta^{-1}\alpha\eta$,
thus $\beta$ is conjugate to $\alpha$.
\end{proof}

Hence, if we want to test whether two braids are conjugate, and we know a $k$th root or the
$k$th power of each one, we just need to test if these roots or powers are conjugate.

This paper is structured as follows. In Section~\ref{S:NT} we give the basic notions and
results from Nielsen-Thurston theory applied to braids. In Section~\ref{S:regular} we study
in more detail a particular case of reducible braids, called {\em reducible braids in
regular form}, that we introduce to simplify the proof of Theorem~\ref{T:main}. This proof
is given in Section~\ref{S:proof}.

\section{Nielsen-Thurston theory}\label{S:NT}

In the same way as isotopy classes of homeomorphisms of surfaces can be classified into
periodic, reducible or pseudo-Anosov \cite{Thclass,FLP}, one has an analogous classification
for braids~\cite{BLM,GMW}.

A braid $\alpha$ is said to be {\em periodic} if it is a root of a power of $\Delta^2$,
where $\Delta=\sigma_1 (\sigma_2\sigma_1)(\sigma_3\sigma_2\sigma_1)\cdots
(\sigma_{n-1}\cdots\sigma_1)$ is Garside's half twist. That is, $\alpha$ is periodic if
$\alpha^k=\Delta^{2m}$ for some nonzero integers $k$ and $m$.

A braid $\alpha$ is said to be {\em reducible} if it preserves (up to isotopy) a family of
disjoint nontrivial simple closed curves on the $n$-punctured disc. Here `nontrivial' means
not isotopic to the boundary but enclosing at least two punctures. Such an invariant family
of curves is called a {\em reduction system}. There exists a {\em canonical reduction
system} $\CRS(\alpha)$ (see~\cite{BLM,Ivanov}), which is the union of all nontrivial curves
$C$ satisfying the following two conditions:
 \begin{enumerate}

 \item $C$ is preserved by some power of $\alpha$.

 \item Any curve $C'$ having nontrivial geometric intersection with $C$ is not preserved by
 any power of $\alpha$.

 \end{enumerate}
It is known that, if $\alpha$ is reducible, then $\CRS(\alpha)=\emptyset$ if and only if
$\alpha$ is  periodic. Hence every reducible, non-periodic braid has a nontrivial canonical
reduction system.

 Finally, a braid $\alpha$ is {\em pseudo-Anosov} if it is neither periodic nor reducible.
In this case \cite{Thclass} there exist two projective measured foliations of the disc,
${\cal F}^u$ and ${\cal F}^s$, which are preserved by $\alpha$. Moreover, the action of
$\alpha$ on ${\cal F}^u$ (the unstable foliation) scales its measure by a real factor
$\lambda >1$, while the action on ${\cal F}^s$ (the stable foliation) scales its measure by
$\lambda^{-1}$. These two foliations and the scaling factor $\lambda$ (called the {\em
stretch factor}), are uniquely determined by $\alpha$.

 Conversely, suppose that a braid $\alpha$ preserves two measured foliations, scaling their
measures by $\lambda$ and $\lambda^{-1}$. One has the following: If $\lambda>1$ then
$\alpha$ is pseudo-Anosov, and if $\lambda=1$ then $\alpha$ is periodic
(see~\cite{Ivanov}).

In order to prove Theorem~\ref{T:main}, we need to show the following results. Although
they are all well-known, we include some short proofs.

\begin{lemma}
 If $\alpha\in B_n$ is periodic, then $\alpha^k$ is periodic, for every $k\neq 0$.
\end{lemma}

\begin{proof}
  There is some $t\neq 0$ such that $\alpha^t$ is a power of $\Delta^2$. Hence
$(\alpha^k)^t=(\alpha^t)^k$ is also a power of $\Delta^2$, thus $\alpha^k$ is periodic.
\end{proof}

\begin{lemma}\label{L:reducible}
If $\alpha\in B_n$ is reducible and not periodic, then $\alpha^k$ is also reducible and
not periodic, for every $k\neq 0$. Moreover, $\CRS(\alpha)=\CRS(\alpha^k)$.
\end{lemma}

\begin{proof} A curve is preserved by a power of $\alpha$ if and only if it is preserved by
a power of $\alpha^k$. Hence, from the definition of the canonical reduction system of a
braid, one has $\CRS(\alpha)=\CRS(\alpha^k)$. Since this family of curves is nonempty, it
also follows that $\alpha^k$ is reducible and not periodic.
\end{proof}

\begin{lemma}\label{L:pA}
If $\alpha\in B_n$ is pseudo-Anosov, with projective foliations ${\cal F}^u$ and ${\cal
F}^s$, and stretch factor $\lambda$, then for every $k\neq 0$, $\alpha^k$ is also
pseudo-Anosov, with projective foliations ${\cal F}^u$ and ${\cal F}^s$, and stretch
factor $\lambda^{|k|}$.
\end{lemma}

\begin{proof}
  This is a straightforward consequence of the definitions.
\end{proof}

\begin{corollary}\label{C:type}
If $\alpha,\beta\in B_n$ are such that $\alpha^k=\beta^k$, then $\alpha$ and $\beta$ are
of the same Nielsen-Thurston type.
\end{corollary}

\begin{proof}
 By the above lemmas, the Nielsen-Thurston type of $\alpha$ (resp. $\beta$) is the same as the type of
 $\alpha^k$ (resp $\beta^k$). Since $\alpha^k=\beta^k$, their types coincide.
\end{proof}

\section{Reducible braids in regular form}\label{S:regular}

The most difficult case in the proof of Theorem~\ref{T:main} occurs when $\alpha$ and
$\beta$ are reducible. Hence, we will study this kind of braid in more detail in this
section. More precisely, we will define a special type of reducible braids, called reducible
braids in {\em regular form}, which are easier to handle if we care about conjugacy. They
were defined in~\cite{GMW} to study centralizers of braid. It is also shown in~\cite{GMW}
that every reducible, non-periodic braid can be conjugated to another one in regular form.
We will repeat that construction here since we need it for our purposes. Later we will give
necessary and sufficient conditions for two braids in regular form to be conjugate. This
will allow us to simplify the proof of Theorem~\ref{T:main}.

\begin{figure}[ht!]
\centerline{\includegraphics{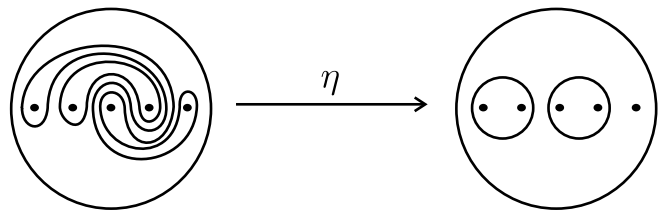}} \caption{Canonical reduction systems can
be simplified} \label{F:simplifyCRS}
\end{figure}

First we will fix a reducible, non-periodic braid $\alpha$. We know that $\CRS(\alpha)$ is
nonempty, but the curves forming $\CRS(\alpha)$ may be rather complicated. If we conjugate
$\alpha$ by some element $\eta$, the canonical reduction system of $\eta^{-1} \alpha \eta$
will be $\eta(\CRS(\alpha))$ (here $\eta$ is considered as an automorphism of the punctured
disc). We can then choose a braid $\eta$ which sends $\CRS(\alpha)$ to the simplest possible
family of closed curves: a family of circles, centered at the real axis (each circle will
enclose more than one and less than $n$ punctures). In other words, up to conjugacy we can
suppose that $\CRS(\alpha)$ is a family of circles (see Figure~\ref{F:simplifyCRS}).

\begin{figure}[ht!]
\centerline{
 \includegraphics{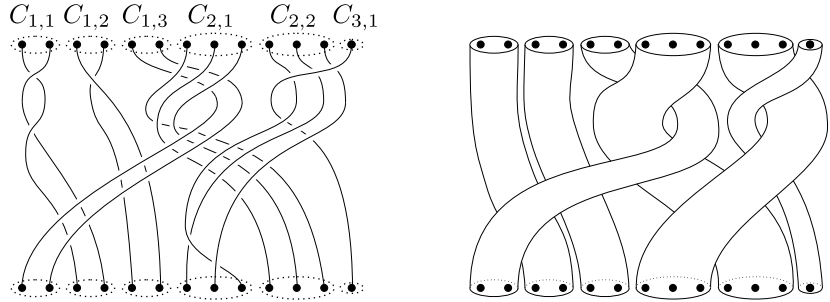} %
} \caption{Example of a reducible braid $\eta$, and its corresponding tubular braid
$\widehat{\eta}$. Notice the indices of the circles of ${\cal C}$.}\label{F:tubular}
\end{figure}

Now we can decompose the punctured disc $D$ in the following way (taken from~\cite{GMW}):
Let ${\cal C}$ be the set of outermost circles of $\CRS(\alpha)$. This set is preserved by
$\alpha$, and we can distinguish the different orbits of circles under $\alpha$. We denote
these orbits by ${\cal C}_1,\ldots, {\cal C}_t$, and the circles forming ${\cal C}_i$ by
$C_{i,1},C_{i,2},\ldots,C_{i,r_i}$. That is, ${\cal
C}=\bigcup_{i=1}^t(\bigcup_{u=1}^{r_i}C_{i,u})$ and $\alpha$ sends $C_{i,u}$ to $C_{i,u+1}$,
where $C_{i,r_i+1}=C_{i,1}$.  In Figure~\ref{F:tubular} we can see an example showing the
notation of these circles. In the examples we will usually number the orbits, and the
circles inside each orbit, from left to right, but this does not need to be true in general:
the only necessary condition is that $\alpha $ sends $C_{i,u}$ to $C_{i,u+1}$. If at some
time we need to stress that these circles, or orbits, belong to $\CRS(\alpha)$, we will
write $C_{i,u}(\alpha)$ or ${\cal C}_i(\alpha )$.

\begin{figure}[ht!]
\centerline{\includegraphics[width=3.5in]{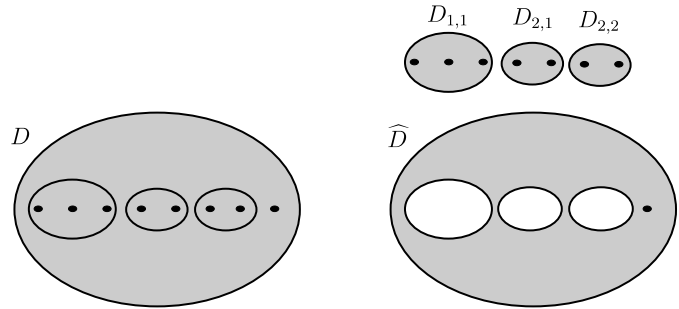}} \caption{Decomposition of the disc along a
canonical reduction system} \label{F:decomposition}
\end{figure}

 Denote by $D_{i,u}=D_{i,u}(\alpha )$ the punctured disc enclosed by the circle $C_{i,u}$, and let
$\widehat{D}=D\backslash(\bigcup_{i,u} {D_{i,u}})$. Notice that $\widehat{D}$ can also be
considered as a punctured disc, if we collapse each hole to a puncture (see
Figure~\ref{F:decomposition}). Hence we have decomposed $D$ into several punctured discs:
$D=\widehat{D}\cup(\bigcup_{i,u}{D_{i,u}})$.

 Now denote by $B_{\cal C}$ the subgroup of $B_n$ formed by those braids that
preserve ${\cal C}$ setwise (maybe permuting the curves that enclose the same number of
punctures). Every braid in $B_{\cal C}$, considered as an automorphism of $D$, induces
automorphisms (braids) on $\widehat{D}$ and on every $D_{i,u}(\alpha)$. More precisely, let
$m$ be the number of punctures in $\widehat{D}$. We can define the map $p:\: B_{\cal C}
\rightarrow B_m$ that sends any $\eta\in B_{\cal C}$ to $\widehat{\eta}$, the braid
corresponding to the automorphism induced by $\eta$ on $\widehat{D}$. It is easy to show
that $p$ is a homomorphism. The braid $p(\eta)=\widehat{\eta}$ is called the {\em tubular}
braid associated to $\eta$.

 In the same way, given $\eta\in B_{\cal C}$, we can define for $i=1,\ldots,t$ and
$u=1,\ldots,r_i$  the braid $\eta_{i,u}=\eta_{i,u,\alpha }$  induced by $\eta$ on the disc
$D_{i,u}(\alpha )$. That is, $\eta_{i,u}$ is the homeomorphism
$\eta_{i,u}:\:D_{i,u}\rightarrow D_{i,u+1}$ induced by $\eta$. These braids are called the
{\em interior} braids of $\eta$. Hence every braid in $B_{\cal C}$ can be decomposed into
one tubular braid and several interior braids (one for each circle in ${\cal C}$).

One can also see this decomposition in a three dimensional picture. If we look at a braid
$\eta\in B_{\cal C}$ in the cylinder $D\times [0,1]$, where the movements of the punctures
are represented as strands, then the movements of the circles $C_{i,u}$ correspond to
`tubes'. If we forget the strands inside the tubes, we obtain the tubular braid
$\widehat{\eta}$, where the solid tubes can be regarded as fat strands. On the other hand,
the strands inside the tube that starts at $C_{i,u}$ and ends at $C_{i,u+1}$, correspond to
the interior braid $\eta_{i,u}$. See an example in Figure~\ref{F:tubular}.

 Notice that every $\eta\in B_{\cal C}$ is completely determined by the braids
 $\widehat{\eta}$ and $\eta_{i,u}$. We can now define a particulary simple kind
 of braid. The definition is taken from~\cite{GMW}:

 \begin{definition}\label{D:regular}
 Let $\alpha$ be a reducible, non-periodic braid, whose canonical reduction system
 is a family of circles. With the above notations, we say that
  $\alpha$ is in {\em regular form} if it satisfies the
 following conditions:
 \begin{enumerate}
    \item For $i=1,\ldots, t$, the interior braids $\alpha_{i,u}$ are all trivial, except
    possibly $\alpha_{i,r_i}$, which is denoted $\alpha_{[i]}$.

    \item For $i,j\in \{1,\ldots,t\}$, the interior braids $\alpha_{[i]}$ and $\alpha_{[j]}$
    are either equal or non-conjugate.

 \end{enumerate}
 \end{definition}

\begin{figure}[ht!]
\centerline{
 \includegraphics{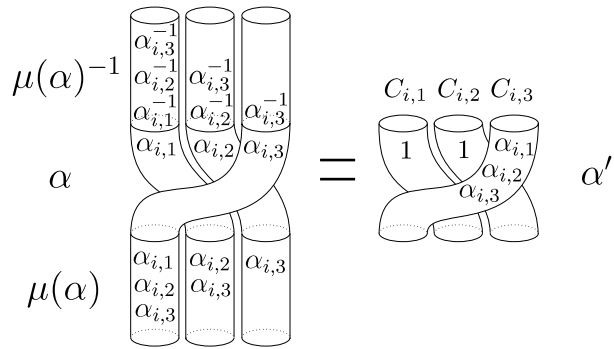} %
} \caption{The conjugation of $\alpha$ by $\mu(\alpha)$ simplifies interior braids.}
\label{F:transfuse}
\end{figure}

It is shown in~\cite{GMW} that every reducible, non-periodic braid $\alpha$ is conjugate to
another one in regular form (although regular forms are not unique, that is, there could be
more than one braid in regular form conjugate to $\alpha$, as we shall see). The precise
conjugation shown in~\cite{GMW} is the following.  We define $\mu(\alpha)\in B_{\cal C}$ as
the braid whose tubular braid $\widehat{\mu(\alpha)}$ is trivial (vertical tubes), and whose
interior braids are the following:
$\mu(\alpha)_{i,u,\alpha}=\alpha_{i,u}\alpha_{i,u+1}\cdots \alpha_{i,r_i}$, for
$i=1,\ldots,t$ and $u=1,\ldots,r_i$. It is an easy exercise to show that
$\alpha'=\mu(\alpha)^{-1} \alpha \mu(\alpha)$ has the same tubular braid as $\alpha$, and
its interior braids are all trivial, except possibly $\alpha'_{i,r_i}=\alpha'_{[i]}=
\alpha_{i,1}\alpha_{i,2}\cdots \alpha_{i,r_i}$. In other words, conjugating by $\mu(\alpha)$
we `transfuse' all interior braids of $\alpha$ in every ${\cal C}_i$ to its last tube (see
Figure~\ref{F:transfuse}). Hence $\alpha'$ satisfies the first condition of
Definition~\ref{D:regular}.

\begin{figure}[ht!]
\centerline{
 \includegraphics{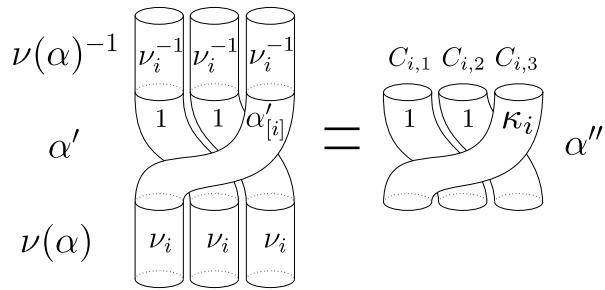} %
} \caption{The conjugation by $\nu(\alpha)$ replaces $\alpha'_{[i]}$ by $\kappa_i$.}
\label{F:nu_alpha}
\end{figure}

Now consider the interior braids $\alpha'_{[1]},\ldots, \alpha'_{[t]}$. For every
$i=1,\ldots,t$, choose one representative $\kappa_i$ of the conjugacy class of
$\alpha'_{[i]}$, in such a way that if $\alpha'_{[i]}$ is conjugate to $\alpha'_{[j]}$, then
$\kappa_i=\kappa_j$. For $i=1,\ldots,t$, choose a braid $\nu_i$ that conjugates
$\alpha'_{[i]}$ to $\kappa_i$. Then we define $\nu(\alpha)\in B_{\cal C}$ as the braid whose
tubular braid $\widehat{\nu(\alpha)}$ is trivial, and whose interior braids are the
following: $\nu(\alpha)_{i,u,\alpha}=\nu_i$ for $i=1,\ldots,t$ and $u=1,\ldots,r_i$. If we
now conjugate $\alpha'$ by $\nu(\alpha)$, then every $\alpha'_{[i]}$ is replaced by
$\kappa_i$ (see Figure~\ref{F:nu_alpha}), that is, $\alpha''=\nu(\alpha)^{-1} \alpha'
\nu(\alpha)$ satisfies the two conditions of Definition~\ref{D:regular}, thus $\alpha''$ is
in regular form.

Therefore every reducible, non-periodic braid $\alpha$ can be conjugated to a braid in
regular form, and one possible conjugating braid is $\mu(\alpha)\nu(\alpha)$. Notice that
$\mu(\alpha)$ depends on our choice of the last tube of each ${\cal C}_i$, and $\nu(\alpha)$
depends on our choice of a representative for the conjugacy class of each $\alpha'_{[i]}$.
Hence, in general, the element $\alpha''$ in regular form conjugate to $\alpha$ is not
unique, though the interior braid $\alpha''_{[i]}=\kappa_i$ is conjugate to
$\alpha_{i,1}\cdots \alpha_{i,r_i}$ in any case.

Next we need a result to test whether two reducible, non-periodic braids in regular form,
having the same canonical reduction system, are conjugate. An obvious necessary condition is
that their tubular braids be conjugate. A necessary and sufficient condition is given by:

\begin{proposition}\label{P:regular_conjugate}
Let $\alpha$ and $\beta$ be two reducible, non-periodic braids in regular form, such that
$\CRS(\alpha)=\CRS(\beta)$ is a family of circles. Then $\alpha$ and $\beta$ are conjugate
if and only if there exists a braid $\eta$ that conjugates $\widehat{\alpha}$ to
$\widehat{\beta}$ such that, if $\eta$ sends the orbit ${\cal C}_i$ of $\widehat{\alpha}$ to
the orbit ${\cal C}_j$ of $\widehat{\beta}$, then $\alpha_{[i]}$ is conjugate to
$\beta_{[j]}$.
\end{proposition}

\begin{proof}
Suppose that $\alpha$ and $\beta$ are conjugate, and let $\xi$ be a conjugating braid, that
is, $\xi^{-1}\alpha \xi=\beta$. It is not difficult to show  (and it is shown
in~\cite{Ivanov}), that in this case $\xi(\CRS(\alpha ))=\CRS(\beta )$. But $\CRS(\beta) =
\CRS(\alpha )$, hence $\xi$ preserves $\CRS(\alpha )$, thus  $\xi\in B_{\cal C}$, where
${\cal C}$ is the set of outermost circles of $\CRS(\alpha)=\CRS(\beta)$. We can then apply
$p$ to all factors of the above equality, and we obtain
$\widehat{\xi}^{-1}\widehat{\alpha}\widehat{\xi}=\widehat{\beta}$, hence $\widehat{\alpha}$
and $\widehat{\beta}$ are conjugate by $\widehat{\xi}$.

Now focus on the permutations induced by $\widehat\alpha $ and $\widehat\beta$ on their base
points $Q_1,\ldots,Q_m$. Consider one cycle $(Q_{i_1},\ldots,Q_{i_r})$ of the permutation
induced by $\widehat\alpha $. Conjugation by $\widehat\xi $ sends it to $(\widehat\xi
(Q_{i_1}),\ldots,\widehat\xi (Q_{i_r}))$, which is a cycle of the permutation induced by
$\beta $ (this is a general property in symmetric groups). Since these cycles correspond to
the orbits of circles of $\alpha $ and $\beta $, then $\xi $ must send any orbit ${\cal
C}_i(\alpha)$ to an orbit ${\cal C}_j(\beta)$.


Consider an orbit ${\cal C}_i(\alpha)$, which is sent to ${\cal C}_j(\beta )$ by conjugation
by $\xi $. The number of circles in each of these two orbits must coincide, so we call it
$r$. Now, since $\xi ^{-1} \alpha \xi =\beta $, one has $\xi^{-1} \alpha^r \xi =\beta^r$.
But $\alpha^r$ is a braid that preserves each circle $C_{i,u}(\alpha )$, and the interior
braid corresponding to any of these circles is $\alpha_{[i]}$ (recall that $\alpha$ is in
regular form). In the same way, $\beta^r$ preserves each circle $C_{j,u}(\beta )$, and the
interior braid corresponding to any of these circles is $\beta_{[j]}$. Choose then some
$C_{i,u}(\alpha)$; it will be sent by $\xi$ to some $C_{j,v}(\beta)$. Then one has:
$$
\beta_{[j]}=(\beta^r)_{j,v,\beta}=\xi_{i,u,\alpha}^{-1}\; (\alpha^{r})_{i,u,\alpha} \;
\xi_{i,u,\alpha} = \xi_{i,u,\alpha}^{-1}\; \alpha_{[i]}\; \xi_{i,u,\alpha}.
$$
Hence $\alpha_{[i]}$ and $\beta_{[j]}$ are conjugate. Therefore, the stated condition is
satisfied by taking $\eta=\widehat{\xi}$.

Conversely, suppose that there exists $\eta$ satisfying the above condition: For every
$i=1,\ldots,t$, the braid $\eta$ sends the orbit ${\cal C}_i$ of $\widehat{\alpha}$ to an
orbit ${\cal C}_j$ of $\widehat{\beta}$, and $\alpha_{[i]}$ is conjugate to $\beta_{[j]}$.
This implies that $\alpha_{[i]}$ and $\beta_{[j]}$ have the same number of strands, that is,
all circles in ${\cal C}_i(\alpha)$, and in ${\cal C}_j(\beta)$, enclose the same number of
punctures. This condition allows us to define a braid $\xi\in B_{\cal C}$ such that
$\widehat{\xi}=\eta$: it suffices to consider the only braid in $B_{\cal C}$ whose tubular
braid is $\eta$ and whose interior braids are all trivial (with the suitable number of
strands).

If we conjugate $\alpha$ by $\xi$, we obtain a braid $\alpha'\in B_{\cal C}$, whose orbits
of circles coincide with those of $\beta$. Moreover, if conjugation by $\xi$ sends ${\cal
C}_i(\alpha)$ to ${\cal C}_j(\beta)$, then in  ${\cal C}_j(\beta)$ there is just one
nontrivial interior braid of $\alpha'$, which is precisely $\alpha_{[i]}$. But this
nontrivial interior braid does not lie, in general, in the last tube of ${\cal C}_j(\beta)$.
Anyway, we know how to conjugate $\alpha'$ to a braid in regular form, with the same
labelling of circles as $\beta$: First, we define the braid $\mu(\alpha')$, such that
$\mu(\alpha')^{-1}\alpha' \mu(\alpha')$ has its nontrivial interior braids (which are still
equal to $\alpha_{[i]}$) in the same tubes as $\beta$. Then, if we take $\beta_{[i]}$ as the
representative for the conjugacy class of $\alpha_{[i]}$, this determines a braid
$\nu(\alpha')$. Conjugating $\alpha'$ by $\mu(\alpha')\nu(\alpha')$, we obtain a braid whose
tubular braid is $\widehat{\beta}$, and whose interior braids coincide with those of
$\beta$, hence we obtain $\beta$. We have then shown that $\alpha$ is conjugate to $\beta$,
and the result follows.
\end{proof}

\section{Proof of Theorem~\ref{T:main}}\label{S:proof}

Suppose that $\alpha,\beta\in B_n$ are such that $\alpha^k=\beta^k$ for some $k\neq 0$. By
Corollary~\ref{C:type}, we can distinguish three cases, depending whether $\alpha$ and
$\beta$ are periodic, reducible or pseudo-Anosov.

\sh{4.1\qua $\alpha$ and $\beta$ are periodic}

\begin{figure}[ht!]
\centerline{\includegraphics{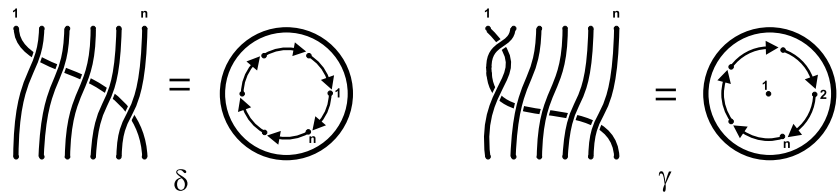}}
\caption{The periodic braids $\delta$ and $\gamma$.}
\label{F:delta_gamma}
\end{figure}

\noindent In this case we shall use a well known result that characterises periodic braids. Consider the
following braids in $B_n$: $\delta=\sigma_1\sigma_2\cdots \sigma_{n-1}$ and
$\gamma=\sigma_1^2\sigma_2\cdots \sigma_{n-1}$. They are drawn in Figure~\ref{F:delta_gamma}.
It is easy to see that $\delta^n=\gamma^{n-1}=\Delta^2$, hence they are both periodic.
The classification result is the following.

\begin{theorem}[de Ker\'ekj\'art\'o~\cite{deK,coko}, Eilenberg~\cite{Eilenberg}] \label{T:kerekjarto}
 Every periodic braid is conjugate either to a power of $\delta$ or to a power
of $\gamma$.
\end{theorem}

Now we will show that $\alpha$ and $\beta$ are conjugate to the same power of $\delta$ or
$\gamma$, so they are conjugate to each other.

If we write every braid in terms of Artin generators, and we notice that the defining
relations in the Artin presentation are homogeneous, it follows that the exponent sum
$s(\eta)$ of a braid $\eta$ is an invariant of its conjugacy class. Notice that
$s(\delta)=n-1$, while $s(\gamma)=n$. Hence, every conjugate of $\delta^t$ has exponent sum
$(n-1)t$, and every conjugate of $\gamma^t$ has exponent sum $nt$. The converse is also
true:

\begin{proposition}
  Let $\alpha$ be a periodic braid. If $s(\alpha)=(n-1)t$ for some $t$, then $\alpha$ is
conjugate to $\delta^t$. If $s(\alpha)=nt$ for some $t$, then $\alpha$ is conjugate to
$\gamma^t$.
\end{proposition}

\begin{proof}
  If $s(\alpha)=(n-1)t$ and $t$ is not a multiple of $n$, then $\delta^t$ is the only power
of $\delta$ or $\gamma$ with the same exponent sum as $\alpha$, hence by
Theorem~\ref{T:kerekjarto}, $\alpha$ is conjugate to $\delta^t$. The same happens if
$s(\alpha)=nt$, and $t$ is not a multiple of $n-1$: in this case $\alpha$ is conjugate to
$\gamma^t$. It remains to show what happens when $s(\alpha)=n(n-1)q$, for some $q$. In this
case $\alpha$ could be conjugate either to $\delta^{nq}$ or to $\gamma^{(n-1)q}$. But
$\delta^{nq}=\gamma^{(n-1)q}=\Delta^{2q}$, hence both statements are true.
\end{proof}

Coming back to our problem, notice that $s(\alpha) k = s(\alpha^k)=s(\beta^k)=s(\beta) k$,
so $s(\alpha)=s(\beta)$, since $k\neq 0$. Hence, by the above proposition, $\alpha$ and
$\beta$ are conjugate to the same power of $\delta$ or $\gamma$, as we wanted to show. The
result is thus true if $\alpha$ and $\beta$ are periodic.

\sh{4.2\qua $\alpha$ and $\beta$ are pseudo-Anosov}

\noindent Let ${\cal F}^u$ and ${\cal F}^s$ be the projective foliations corresponding to
$\alpha$, and let $\lambda$ be its stretch factor. By Lemma~\ref{L:pA}, the foliations
and stretch factor corresponding to $\alpha^k=\beta^k$ are ${\cal F}^u$, ${\cal F}^s$ and
$\lambda^{|k|}$. Therefore, those corresponding to $\beta$ are ${\cal F}^u$, ${\cal F}^s$
and $\lambda$.

 We will show that, in this case, $\alpha$ and $\beta$ commute: their commutator,
$\rho=\alpha \beta \alpha^{-1} \beta^{-1}$ preserves ${\cal F}^u$ and ${\cal F}^s$, and its
stretch factor is 1. Therefore, $\rho$ is periodic, thus conjugate to a power of $\delta$ or
$\gamma$. Since the exponent sum $s(\rho)=0$ then $\rho=1$, so $\alpha$ and $\beta$ commute.

 Now since $\alpha^k=\beta^k$ one has $\alpha^k \beta^{-k}=1$. But since $\alpha$ and
$\beta$ commute, $\alpha^k \beta^{-k}=(\alpha \beta^{-1})^k=1$. This implies that $\alpha
\beta^{-1}$ is a torsion element of $B_n$, but since $B_n$ is torsion free, it follows that
$\alpha \beta^{-1}=1$ hence $\alpha=\beta$. Therefore, if $\alpha$ and $\beta$ are
pseudo-Anosov, then not only $\alpha$ is conjugate to $\beta$: they coincide.

\sh{4.3\qua $\alpha$ and $\beta$ are reducible, not periodic}

\noindent We will show this case by induction on the number of strands. Since the case
$n=2$ has already been studied (all braids on two strands are periodic), we can suppose
that $n>2$, and that our result is true for every pair of braids with less than $n$ strands.

 Since $\alpha$ and $\beta$ are reducible and not periodic, we know that their canonical
reduction systems are non-empty and, by Lemma~\ref{L:reducible}, they must coincide since
they are both equal to $\CRS(\alpha^k)=\CRS(\beta^k)$.

We will now conjugate $\alpha$ and $\beta$ by some braid $\eta$, in order to simplify their
canonical reduction system. This can be done due to the following fact: since
$\alpha^k=\beta^k$, we have that $(\eta^{-1}\alpha \eta)^k= \eta^{-1}\alpha^k \eta =
\eta^{-1}\beta^k \eta=(\eta^{-1}\beta \eta)^k$. If we then show that $\eta^{-1}\alpha \eta$
and $\eta^{-1}\beta \eta$ are conjugate, then $\alpha$ and $\beta$ will also be conjugate.
Therefore, it suffices to show the result for suitable conjugates of $\alpha$ and $\beta$,
hence we can suppose that $\CRS(\alpha)=\CRS(\beta)$ is a family of circles.

 As usual, we denote by ${\cal C}$ the set of outermost circles of $\CRS(\alpha)=\CRS(\beta)$.
Since $\alpha,\beta\in B_{\cal C}$, we can study their tubular and interior braids. First,
since $\alpha^k=\beta^k$ one has $\widehat{\alpha}^k=\widehat{\beta}^k$. Moreover, since
${\cal C}$ is formed by the outermost circles of $\CRS(\alpha)=\CRS(\beta)$, this implies
that $\CRS(\widehat{\alpha})=\CRS(\widehat{\beta})=\emptyset$. That is, $\widehat{\alpha}$
and $\widehat{\beta}$ are either pseudo-Anosov or periodic. We will treat these two cases
separately.

\rk{\bf $\widehat{\alpha}$ and $\widehat{\beta}$ are pseudo-Anosov} In
this case, since $\widehat{\alpha}^k=\widehat{\beta}^k$, we have already shown that
$\widehat{\alpha}=\widehat{\beta}$. Hence we can label the circles of ${\cal C}$ in the same
way for both braids: ${\cal C}=\bigcup_{i=1}^t({\cal
C}_{i})=\bigcup_{i=1}^t(\bigcup_{u=1}^{r_i}(C_{i,u}))$. We will show that $\alpha$ and
$\beta$ are conjugate by conjugating them to the same braid in regular form.

Recall that if we conjugate $\alpha$ by $\mu(\alpha)\nu(\alpha)$, we obtain a braid
$\alpha''$, in regular form, such that $\alpha''_{[i]}$ (the interior braid starting at
$C_{i,r_i}$) is conjugate to $(\alpha_{i,1}\alpha_{i,2}\cdots\alpha_{i,r_i})$ for
$i=1,\ldots,t$. In the same way, if we conjugate $\beta$ by $\mu(\beta)\nu(\beta)$, we
obtain a braid $\beta''$, in regular form, such that $\beta''_{[i]}$ is conjugate to
$(\beta_{i,1}\beta_{i,2}\cdots\beta_{i,r_i})$ for $i=1,\ldots,t$. Therefore, $\alpha''$ and
$\beta''$ are two braids in regular form, whose tubular braids coincide
($\widehat{\alpha''}=\widehat{\alpha}=\widehat{\beta}=\widehat{\beta''}$), and whose
nontrivial interior braids are placed into the same tubes. It just remains to show that we
can take $\alpha''_{[i]}=\beta''_{[i]}$. But $\alpha''_{[i]}$ and $\beta''_{[i]}$ are just
representatives of their conjugacy classes, hence it suffices to show that
$(\alpha_{i,1}\alpha_{i,2}\cdots\alpha_{i,r_i})$ is conjugate to
$(\beta_{i,1}\beta_{i,2}\cdots\beta_{i,r_i})$ for $i=1,\ldots,t$. In order to prove this, we
can forget about $\alpha''$ and $\beta''$, and look at $\alpha$ and $\beta$, and their
tubular and interior braids, as follows.

 We know that $\widehat{\alpha}^k=\widehat{\beta}^k$. Up to raising this braid to a
suitable power, we can suppose that $\widehat{\alpha}^k$ is a pure braid (its corresponding
permutation is trivial). Hence, for $i=1,\ldots,t$, the length $r_i$ of the orbit ${\cal
C}_i$ must be a divisor of $k$, say $r_i p_i=k$. We will now look at the interior braids of
$\alpha^k$ and $\beta^k$. If we raise $\alpha$ to the power $r_i$, then the interior braid
$(\alpha^{r_i})_{i,1,\alpha }=\alpha_{i,1}\cdots \alpha_{i,r_i}$. Hence, if we further raise
it to the power $p_i$, we obtain $(\alpha^k)_{i,1,\alpha }=(\alpha_{i,1}\cdots
\alpha_{i,r_i})^{p_i}$. In the same way, one has $(\beta^k)_{i,1,\beta }=(\beta_{i,1}\cdots
\beta_{i,r_i})^{p_i}$. Since $\alpha^k=\beta^k$, and $C_{i,1}(\alpha)=C_{i,1}(\beta)$, it
follows that $(\alpha_{i,1}\cdots \alpha_{i,r_i})^{p_i}=(\beta_{i,1}\cdots
\beta_{i,r_i})^{p_i}$, which yields, by the induction hypothesis, that
$(\alpha_{i,1}\cdots\alpha_{i,r_i})$ is conjugate to $(\beta_{i,1}\cdots\beta_{i,r_i})$, as
we wanted to show.

\rk{$\widehat{\alpha}$ and $\widehat{\beta}$ are periodic} This
time $\widehat{\alpha}^k=\widehat{\beta}^k$ implies that $\widehat{\alpha}$ and
$\widehat{\beta}$ are conjugate (we have shown this for periodic braids), but not
necessarily equal. Since they are periodic and conjugate, then they are both conjugate to
the same power of $\delta$ or $\gamma$. We can also assume -- as above -- that
$\widehat{\alpha}^k=\widehat{\beta}^k$ is a pure braid, that is, a power of $\Delta^2$ (with
$m$ strands).

Suppose that $\widehat{\alpha}$ and $\widehat{\beta}$ are conjugate to some power of
$\delta$. In this case, all the orbits of outermost circles of $\alpha$ and $\beta$ (i.e.
orbits of points of $\widehat{\alpha}$ and $\widehat{\beta}$) have the same length, say $r$.
And $k$ is a multiple of $r$, say $k=pr$. But the orbits of $\alpha$ and $\beta$ (resp.
$\widehat{\alpha}$ and $\widehat{\beta}$) do not necessarily coincide.

We will first conjugate $\alpha$ to a braid in regular form. Let us choose, from now on, a
representative for each conjugacy class of braids. Then, for $i=1,\ldots,t$, let $\kappa_i$
be the representative of the conjugacy class of $(\alpha_{i,1}\cdots\alpha_{i,r})$. Notice
that, by construction, if $\kappa_i$ and $\kappa_j$ are conjugate, then $\kappa_i=\kappa_j$.
Recall that there exists a braid $\mu(\alpha)\nu(\alpha)$ which conjugates $\alpha $ to
$\alpha''$, in regular form, whose nontrivial interior braids are
$\kappa_1,\ldots,\kappa_t$.

We will now see that the list $\kappa_1,\ldots,\kappa_t$ is completely determined by the
interior braids of $\alpha^k$. Indeed, take some circle $C_{i,u}(\alpha)$. Since the orbit
${\cal C}_i(\alpha)$ has length $r$, and $k=pr$, it follows that
$(\alpha^k)_{i,u,\alpha}=(\alpha_{i,u}\cdots \alpha_{i,r}\alpha_{i,1}\cdots
\alpha_{i,u-1})^p$. Conjugating this braid by $\alpha_{i,u}\cdots \alpha_{i,r}$, we obtain
$(\alpha_{i,1}\cdots\alpha_{i,r})^p$, which is conjugate to $\kappa_i^p$. That is,
$(\alpha^k)_{i,u,\alpha}$ is conjugate to $\kappa_i^p$. Hence, the interior braids of
$\alpha^k$ inside the circles $C_{i,1}(\alpha ),\ldots, C_{i,r}(\alpha)$ are all conjugate
to $\kappa_i^p$. Suppose that they were also conjugate to $\kappa_j^p$, for $j\neq i$. Then
$\kappa_i^p$ would be conjugate to $\kappa_j^p$. But these braids have less than $n$
strands, so by Corollary~\ref{C:main} (that we are allowed to use by the induction
hypothesis), $\kappa_i$ would be conjugate to $\kappa_j$, and then $\kappa_i=\kappa_j$.

In other words, the list of representatives $\kappa_1,\ldots,\kappa_t$, counting
repetitions, is completely determined by the interior braids of $\alpha^k$, as follows.
First we define a partition of ${\cal C}$ by the following property: two circles belong to
the same coset if and only if the corresponding interior braids of $\alpha^k$ are conjugate.
Then the number of circles in every coset is always a multiple of $r$, and the $p$-th roots
of these interior braids determine the representatives $\kappa_1,\ldots,\kappa_t$. An
element $\kappa_i$ appears repeated $q$ times in the list if and only if its corresponding
coset has $qr$ circles.

Now we can repeat the whole construction with $\beta$. We will conjugate $\beta $ to
$\beta''$, in regular form, whose list of nontrivial interior braids is completely
determined by $\beta^k$. But $\alpha^k=\beta^k$, hence the list of nontrivial interior
braids of $\beta''$ is exactly $\kappa_1,\ldots,\kappa_t$.

We then have two braids $\alpha''$ and $\beta''$ in regular form, whose lists of nontrivial
interior braids coincide, and whose tubular braids, $\widehat{\alpha''}=\widehat{\alpha}$
and $\widehat{\beta''}=\widehat{\beta}$, are both conjugate to the same power of $\delta$.
We must show that, in this case, $\alpha''$ and $\beta''$ are conjugate. Up to conjugating
$\alpha''$ and $\beta''$ by elements in $B_{\cal C}$ with trivial interior braids, we can
suppose that $\widehat{\alpha''}= \widehat{\beta''}=\delta^s$, for some $s$. Hence, the
orbits of $\alpha''$ and $\beta''$ coincide, although the corresponding interior braids
could lie in different tubes (even in different orbits).

We can now apply Proposition~\ref{P:regular_conjugate}. $\alpha''$ and $\beta''$ will be
conjugate if it exists a braid $\eta$ that conjugates $\widehat{\alpha''}$ to
$\widehat{\beta''}$ (that is, a braid that commutes with $\delta^s$) such that, if
$\eta({\cal C}_i)={\cal C}_j$, then $\kappa_i=\kappa_j$. In other words, we need an element
of the centralizer of $\delta^s$ which permutes the orbits of $\delta^s$ in the appropriate
way. Fortunately, such an element always exists.

The centralizer of a power of $\delta$ has been described in~\cite{BDM} (see
also~\cite{GMW}): consider the braid $\delta^s=(\sigma_1\cdots\sigma_{m-1})^s$ on $m$
strands. We will now denote by ${\cal C}_i$ the orbit (of points) induced by $\delta^s$ that
starts at the point $P_i$. That is, ${\cal C}_i=
\{P_i,P_{s+i},P_{2s+i}\ldots,P_{(r-1)s+i}\}$, where the indices are taken modulo $m$. Then,
for $i=1,\ldots,t-1$, consider the braid $S_i=\left(\sigma_i \delta^s\right)^r$. This braid
commutes with $\delta^s$ and permutes the orbits ${\cal C}_i$ and ${\cal C}_{i+1}$.
Therefore, taking products of the elements $S_i$, we can obtain any desired permutation of
the orbits. Hence the braid $\eta$ required by Proposition~\ref{P:regular_conjugate} exists,
so $\alpha''$ is conjugate to $\beta''$. Therefore, if $\widehat{\alpha}$ and
$\widehat{\beta}$ are conjugate to a power of $\delta$, then  $\alpha$ is conjugate to
$\beta$, as we wanted to show.

 It remains the case when $\widehat{\alpha}$ and $\widehat{\beta}$ are conjugate to $\gamma^s$, for some
 $s$. Recall that $\widehat{\alpha}$ and $\widehat{\beta}$ have $m$ strands.
 We can suppose that $s$ is not a multiple of $m-1$, since in that case $\gamma^s$ is a power of $\Delta^2$,
 thus a power of $\delta$, and this case has already been studied. Therefore, the orbits of
 $\widehat{\alpha}$ and $\widehat{\beta}$ are as follows: there is one orbit ${\cal C}_1$ of length one,
 and all the other orbits ${\cal C}_2,\ldots,{\cal C}_t$ have the same length, say $r$, where $r>1$.

Now we can apply the same methods as before. First, we can suppose that $\widehat{\alpha}^k=
\widehat{\beta}^k$ is a pure braid, and we denote $p=k/r$. Then the interior braids of
$\alpha^k$ are as follows: one of them equals $(\alpha_{1,1})^k$, and the others are
conjugate to $(\alpha_{i,1}\cdots\alpha_{i,r})^p$ for some $i$. We then choose a
representative $\kappa_1$ for the conjugacy class of $\alpha_{1,1}$, and a representative
$\kappa_i$ for the conjugacy class of $(\alpha_{i,1},\ldots,\alpha_{i,r})$, for $i>1$. This
time, the list $\kappa_1,\ldots,\kappa_t$ is determined by $\alpha^k$ as follows. First we
define a partition of ${\cal C}$ by the following property: two circles of ${\cal C}$ belong
to the same coset if and only if the  interior braids of $\alpha^k$ in their corresponding
tubes are conjugate. Then there is just one coset whose size is not a multiple of $r$, but
congruent to 1 modulo $r$. Take any interior braid of $\alpha^k$ in that coset. Its $k$th
root is conjugate to $\kappa_1$. The $p$th roots of the remaining interior braids of
$\alpha^k$ yield the elements $\kappa_2,\ldots,\kappa_t$, as in the previous case.

 Now we conjugate $\alpha$ to $\alpha''$, in regular form, whose nontrivial interior braids are
 $\kappa_1,\ldots,\kappa_t$. Then we conjugate $\beta$ to $\beta''$, in regular form, whose nontrivial
 interior braids are also $\kappa_1,\ldots,\kappa_t$, since they are determined by $\beta^k=\alpha^k$.
 Notice that, in both cases, $\kappa_1$ is the interior braid of the only orbit of length one. We can now
 conjugate $\alpha''$ and $\beta''$ by suitable braids in $B_{\cal C}$, with trivial interior braids,
 thus we can suppose that $\widehat{\alpha''}=\widehat{\beta''}=\gamma^s$. This time, the interior braid
 $\kappa_1$ is already in the same position (the first tube of $\gamma^s$) for both braids, but the
 remaining interior braids could be in different positions. Applying Proposition~\ref{P:regular_conjugate}
 again, we need a braid $\eta$ that commutes with $\gamma^s$, permuting the orbits ${\cal C}_2,
 \ldots,{\cal C}_t$ as desired. This is always possible, as we are about to see.

  The centralizer of $\gamma^s=(\sigma_1^2\sigma_2\cdots\sigma_{m-1})^s$ has also been studied
 in~\cite{BDM}. For every $i=2,\ldots,t-1$, the braid $T_i=\left(\sigma_i \gamma^s\right)^r$
 commutes with $\gamma^s$ and permutes the orbits ${\cal C}_i$ and ${\cal C}_{i+1}$. Therefore we can
 always define $\eta$ as a product of the $T_i$'s, hence $\alpha''$ and $\beta''$ are conjugate. It
 follows that $\alpha$ and $\beta$ are conjugate, and the theorem is proved.
\medskip

\noindent {\bf Acknowledgements}\qua I am very grateful to Bert Wiest for useful conversations,
and for his comments and suggestions on an early version of this paper. I also want to thank
Tara Brendle, who pointed out some inconsistent notation, and helped to clarify some obscure
parts of the proofs.

The author is partially supported by MCYT, BFM2001-3207 and FEDER.

\Addresses\recd

\end{document}